\documentclass[10]{amsart}
\usepackage[cp1251]{inputenc}
\usepackage{amsfonts}
\usepackage{amssymb,amsmath}

\usepackage[english]{babel} \theoremstyle{definition}

\theoremstyle{remark}

\begin{document}

\title[]{Occasionally attracting compact sets and
compact-supercyclicity}
\author{K. V. Storozhuk}
\address{Sobolev Institute of Mathematics, Akademik
Koptyug Prospect 4, 630090, Novosibirsk, Russia}
%\thanks{Author was supported by the Program "Universites of Russia", Grant 8311}

\date{March, 2006}

%\keywords{$C_0$-semigroup, power bounded operator, supercyclicity,
%almost-periodic representation}

\begin{abstract} If $(T_t)_{t\ge 0}$ is a bounded
$C_0$-semigroup in a Banach space $X$ and there exists a compact
subset $K\subseteq X$ such that
$$
  \liminf_{t\to\infty}\rho(T_tx,K)=0 \ \ \ (\forall x\in X, \ \|x\|\leq 1), $$
then there exists a finite-dimensional subspace $L\subseteq X$ such
that
$$
  \lim_{t\to\infty}\rho(T_tx,L)=0 \ \ \ (\forall x\in X).
$$

If $T:X\to X$ ($X$ is real or complex) is supercyclic and
$(\|T^n\|)_n$ is bounded then $(T^nx)_n$ vanishes for every $x\in
X$.

We define the "compact-supercyclicity". If $\dim X=\infty $ then $X$
has no compact-supercyclic isometries.
\end{abstract}

\maketitle \centerline{\bf Introduction}\vspace{1mm}

Let $X$ be a Banach space, by $B_X$ we denote the unit ball in $X$.
For a subset $Y\subseteq X$ and $x\in X$ we denote by $\rho(x,Y)$
the distance between $x$ and $Y$.

\vspace{2mm} {\bf Definition 1} \ Let $T:X\to X$ be a map. A subset
$K\subseteq X$ is called an {\it attractor} for $T$ if
$$
  \forall x\in B_X\ \lim_{n\to\infty}\rho(T^nx,K)=0.
  \eqno(1)
$$
\vspace{2mm}

The definition of an attractor for a $C_0$-semigroup is similar. It
is known that for a linear power bounded operator (and for a bounded
$C_0$-semigroup) the existence of a compact attractor implies the
existence of an invariant finite-dimensional subspace $L\subseteq X$
and an invariant subspace $X_0\subseteq X$ such that $X=X_0\oplus L$
and the semigroup $(T^n)_n$ is isomorphic to the direct product of
semigroups
$$
  T^n=(T|_{X_0})^n\oplus (T|_L)^n:X_0\oplus L\to X_0\oplus L, \  \forall x_0\in X_0\ T^nx_0\to 0.
  \eqno(2)
$$
This theorem was proved in \cite{LLY} for the Markov semigroups in
$L_1$. Its general case was proved by  Vu \cite{Ph} and Sine
\cite{S}. We call this result the Vu -- Sine theorem.

It turns out that the conclusion of the Vu -- Sine theorem remains
true if there exists only ``occasionally attracting" compact set
$K$:
$$
  \liminf_{t\to\infty}\rho(T_tx,K)=0 \ \ \ (\forall x\in B_X).
  \eqno(1')
$$

The papers \cite{Ph} and \cite{S} use the results of Jacobs \cite{J}
and de Leeuw and Glicksberg \cite{LG} on spectral decomposition of
weakly almost periodic semigroup. We base on a more elementary fact,
the non-emptiness of the essential spectrum.

In the first part of the paper we prove the following theorem:
% hat
%if an occasionally attracting compact set exists for an  \it
%isometry \rm\ of $X$ then $X$ has finite dimension:

\vspace{2mm} {\bf Theorem 1}. \it Let $\dim X=\infty$. For any
isometry $T:X\to X$ there are no occasionally attracting compact
sets.\rm\vskip1mm

The second part of the paper is devoted to the application of
Theorem 1 to the above-mentioned strengthening of the results of
\cite{LLY},\cite{Ph}, \cite{S}.

\vspace{2mm} {\bf Definition 2}. Let $x\in X$ or $x\subseteq X$.
Denote by $O(x)=\cup_{n=0}^\infty T^nx$ the orbit of $x$.

\vspace{2mm} {\bf Definition 3}.  A vector $a$ is called a {\it
returning   vector} if $ \liminf_{n\to
\infty}\|T^na-a\|=0.$\vskip1mm

It is easy to see that if $T$ is power bounded then $a\in X$ is a
returning   vector if and only if $a$ is a limit point of the orbit
of some $x\in X$. \vspace{1mm}

 {\bf Lemma 1} \ {\it Let $T:X\to X$ and
$\|T\|\leq 1$. If $a$ is a returning   vector, then the subspace
$L(a)=\mbox{cl}(\mbox{span}(O(a))$ consists of returning vectors and
$T:L(a)\to L(a)$ is an isometry.} \vspace{1mm}

%Then, by using of Theorem 1, we prove the above-mentioned
%generalizations of the Vu -- Sine theorem and the existence
%of representation~(2).

Then we prove the generalizations of the Vu -- Sine theorem from
theorem 1:

\vspace{1mm} {\bf Theorem 2} \ {\it Let $T:X\to X$ be a power
bounded operator. If there exists a compact set $K$ such that $(1')$
holds, then the semigroup $(T^n)_{n=0}^{\infty}$ is ``asymptotically
finite-dimensional", i.e. there exists an invariant subspace
$L\subseteq X$, $\dim(L)<\infty$ such that, for every $x\in X$,
$\lim_{n\to\infty}\rho(T^nx,L)=0$ and decomposition $(2)$ holds. The
space  $L$ is generated by all returning
  vectors of\ \ $T$.} \vspace{2mm}

{\bf Theorem 3} \ {\it Let $(T_t)_{t\ge 0}$ be a bounded
$C_0$-semigroup in a Banach space $X$. If there exists an
occasionally attracting compact set $K\subseteq X$, then the
semigroup $T$ is asymptotically finite-dimensional.} \vspace{2mm}

The last part of the paper is devoted to another application of
Theorem 1.\vspace{2mm}
%to the theory
%of supercyclisity in real and complex Banach spaces.

Let $X$ be a real or complex infinite-dimensional Banach spaces and
$F\in\{\Bbb R,\Bbb C\}$. An operator $T:X\to X$ is called \it
 supercyclic \rm if there exists a vector $k\in X$ such that the set
$F\cdot O(k)$
% =\{\lambda x\mid \lambda\in F,\ x\in O(k)\}$
is dense in $X$. The corresponding vector $k$ is called \it
supercyclic\rm.

The following results were proved for complex $X$ in \cite{AB} and
\cite{Mi}:

{\bf Theorem 4 \it
%(in the case $F=\Bbb C$ was proved in \cite{AB}).
If \ $T:X\to X$ is isometry, then $T$ is not supercyclic. Moreover,
if $T$ is power bounded and supercyclic, then $T^nx$ vanishes for
every $x\in X$\rm.\vskip1mm

Both \cite{AB} and \cite{Mi} make use of the Godement theorem
\cite{Go}: every isometry of complex $X$ has an invariant proper
closed subspace.

We deduce Theorem 4 (in the real and complex cases) from Theorem 1.
The proof is based on the following lemma:

{\bf Lemma 4}. \it Let $\|T\|\leq 1$. If $T^na\not{\!\!\to}0$ and
there exist $\lambda_k$ and $n_k$ such that $\lambda_kT^{n_k}a\to a$
$($in particular, if $a$ is supercyclic$)$, then $a$ is a returning
vector\rm.

\it Remark\rm. %Through the rest of the introduction we
In \cite{Mi} Miller proved that an isometry of a complex $X$ cannot
even be \it finite-supercyclic\rm\, i.e., for any finite set
$K\subseteq Z$, the set $F\cdot O(K)$ is not everywhere dense in
$Z$. But after that, Peris \cite{Pe} showed that, for locally convex
spaces, finite-supercyclicity is equivalent to supercyclicity. A
weaker property is  \it N-supercyclicity\rm\ \cite{Fe, BoFeSha}. An
operator  $T$ is N-supercyclic if there exists a finite-dimensional
unit ball $B_L\subseteq X$ such that $\Bbb C\cdot O(B_L)$ is dense
in $X$. Following this tradition, we may call $T:X\to X$ {\it
compact-supercyclic} if there is a compact set $K\subseteq X$ such
that $F\cdot O(K)$ is dense in $X$. If $\dim X<\infty$ then each $T$
is compact-supercyclic.

For example, if $T:X\to X$ is an isometry then the presence of
an~occasionally attracting compact set for $T$ is equivalent to
compact-supercyclicity of $T^{-1}$ (cf. the proof of Theorem 4).
Therefore, we can reformulate Theorem 1 as follows:

\it If $\dim X=\infty $ then $X$ has no compact-supercyclic
isometries\rm.

\maketitle \centerline{1 \ \ \bf Proof of Theorem 1}\vspace{1mm}

First we consider the case of a complex $X$. Let $\sigma_{ess}(T)$
be the essential spectrum.
%of the operator $T$
If $\lambda\in\sigma_{ess}(T)$ then $\dim\ker(T-\lambda)=\infty$ or
the $\text{Im}(T-\lambda)$ is not closed in $X$.

We say that a bounded sequence $z_n\in X$ is  {\it sparse} if it
contains no converging subsequence. Let us show that it is possible
to assign to each $\lambda\in\sigma_{ess}(T)$ a sparse sequence of
``approximate eigenvectors"\  $z_n\in B_X$, i.e. such that
$Tz_n-\lambda z_n\to 0$. Borrowing the terminology from the theory
of self-adjoint operators in Hilbert spaces, we call such a sequence
$z_n$ a {\it Weyl sequence}.

\vspace{1mm} {\bf Lemma 2} \ {\it For each
$\lambda\in\sigma_{ess}(T)$, there exists a Weyl sequence $z_n$.}
\vspace{1mm}

{\it Proof}: \ Put $S=(T-\lambda):X\to X$. We have: either $\dim\ker
S=\infty$ or $S(X)$ is not closed in $X$. If $\dim\ker S=\infty$,
then the statement is obvious.

If $\dim\ker S<\infty$ then $\ker S$ has a closed complement
$V\subseteq X$. Consider the operator $S|_V:V\to X$. The kernel of
$S|_V$ is zero and its image $S|_V(V)=S(X)$ is not closed in $X$.
Therefore the inverse operator $(S|_V)^{-1}:S(X)\to V$ is unbounded
and there exists a sequence $z_n\in V$, $\|z_n\|=1$ such that
$Sz_n\to 0$. The sequence $z_n$ has no limit points, since they
would be nonzero elements of the kernel of $S|_V$.$\Box$

\vspace{2mm} If $z_n\in X$ is a sparse sequence and $Tz_n-\lambda
z_n\to 0$ then
$$
  \forall k\in\Bbb N \ \ \|T^k z_n-\lambda T^{k-1}z_n\|=\|Tz_n-\lambda z_n\|\to 0.
  \eqno(3)
$$

Suppose that $K$ is an~occasionally attracting compact set. For each
$n\in\Bbb N$, there exist a number $k_n$ and $a_n\in K$ such that
$\|T^{k_n}z_n-a_n\|<\frac1n.$ Switching to a subsequence, one can
assume that $a_n\to a$ and $\|T^{k_n}z_n-a\|\to 0$, i.e.
$T^{k_n}z_n\to a$. It follows from (3) that $Ta=\lambda a$. In
particular, the $\Bbb Z$-orbit $\{T^na\mid n\in \Bbb Z\}$ of $a$
lies in some one-dimensional subspace $L(a)\subseteq X$. But
$$
  \rho(z_n, L(a))\leq\|z_n-T^{-k_n}a\|=\|T^{k_n}z_n-a\|\to 0, $$
i.e., the sequence $z_n$ approaches a~one-dimensional subspace and
thus cannot be sparse. The theorem is proved in the complex case.

\vspace{2mm} The real case requires the following auxiliary lemma.

\vspace{2mm} {\bf Lemma 3} (an analog of spectrum in real space) \
{\it Let $X$ be a real space and let $T:X\to X$ be a bounded
operator. There exist two numbers $r,s\in \Bbb R$ such that the
operator $S:=T^2+rT+s$ is not bijective. Moreover, if $\dim
X=\infty$, then there exist $r,s\in\Bbb R$ and a sparse sequence
$x_n$ such that $T^2x_n+rTx_n+sx_n\to 0$.} \vspace{2mm}

{\it Proof}: \ Any complex $\lambda$ is a root of the real
polynomial
$$
  S_\lambda(t)=(t-\lambda)(t-\overline \lambda)=t^2-t(\lambda+\overline\lambda)+|\lambda|^2.
$$
Consider the complexification: $T_\Bbb C:X_\Bbb C\to X_\Bbb C$,
$T_\Bbb C(x+iy)=Tx+iTy$. If $\lambda\in\sigma(T_\Bbb C)$, then the
operator $T_\Bbb C-\lambda$ is not bijective, hence the  operator
$S_\lambda(T_{\Bbb C})$ is not bijective either. On the other hand,
the coefficients of the polynomial $S_\lambda$ are real, so
$S_\lambda(T_{\Bbb C})= (S_\lambda(T))_{\Bbb C}$; therefore
$S_\lambda (T):X\to X$ is not bijective as well.

If $\dim X=\infty$, let $\lambda\in\sigma_{ess}(T_\Bbb C)$ and
$z_n=x_n+iy_n\in X_\Bbb C$ be the corresponding Weyl sequence. Then
$S_\lambda(T_{\Bbb C})z_n\to 0$. But then  $S_\lambda(T)x_n\to 0$
and $S_\lambda(T)y_n\to 0$. The sequences $x_n$ and $y_n$ do not
have to be sparse, but if in the sequence $y_n\in X$ of the
imaginary parts of $z_n\in X_{\Bbb C}$ there can be found a
converging subsequence $y_{n_k}\in X$ then the corresponding
subsequence of the real parts $x_{n_k}\in X$ is certainly sparse.
The lemma is proved. $\Box$

{\bf Example}.  {\it If\ $T:\Bbb R^2\to\Bbb R^2$ is the rotation on
$\alpha\in(0,\pi)$, then $T^2-\frac{\sin
2\alpha}{\sin\alpha}T+1=0$.}

Now we finish the proof of Theorem 1 in the real case. Let $x_n$ be
a sparse sequence such that $T^2x_n+rTx_n+sx_n\to 0$. By arguments
as in the proof of the complex case we find  $a\in K$ such that
$T^2a+rTa+sa=0$. The orbit of the vector $a$ belongs to the
two-dimensional subspace $L(a)$, attracting some subsequence in
$x_n$. This contradicts the $x_n$ being sparse. Theorem 1 is
completely proved. \vspace{2mm}

\maketitle \centerline{2\ \ \bf Proofs of  Theorems 2 and
3}\vspace{1mm}

{\it Proof of Lemma 1}: \ Notice that
$\|a\|=\|Ta\|=\|T^2a\|=\ldots$. Indeed, this sequence is
non-increasing. The vector $a$ is returning, therefore this sequence
cannot decrease either. Now, for each $n\in \Bbb N$ the vector
$T^na$ is also returning, such are also linear combinations of these
vectors. So, we have $\|T(x)\|=\|x\|$ for every $x\in L(a)$.
Finally, the set $T(L(a))$ is dense in $L(a)$ and therefore
$T(L(a))=L(a)$. $\Box$

\vspace{1mm} {\it Proof of Theorem 2}: \ Assume that  $\|T\|\leq 1,$
rescaling $X$ by the equivalent norm
$$\|x\|:=\sup\{\|x\|,\|Tx\|,\|T^2x\|,\ldots\}.\eqno(4)$$

For each $x\in B_X$, there is $a\in K$ such that $a$ is the limit
points of the orbit $O(x)$. It is clear that $T^nx\to O(a)$ and $a$
is returning vector.

According to Lemma 1 and Theorem 1, the set $O(a)$
 lies in an invariant finite-dimensional subspace $L(a)$.

%Each $x\in B_X$ is attracted by the orbit $O(a)$ of some
%returning
%  vector $a\in K$. Indeed, let  $a\in K$ be a limit point of the
%orbit of the vector $x$. The vector $a$ is returning. If
%$\|T^nx-a\|<\varepsilon$, then for each $k>0$
%$\|T^{n+k}x-T^ka\|=\|T^k(T^nx-a)\|<\varepsilon$ and $T^na\to O(a)$,
%since $\liminf_{n\to\infty}\|T^nx-a\|=0$.

Different vectors $x_i$ of the unit ball in $X$ can be attracted,
generally speaking, by the orbits of different vectors $a_i\in K$.
It remains to prove that the orbits of all returning   vectors lie
in one and the same finite-dimensional space $L=\oplus L(a_i)$.

Let $a_1,\ldots a_n$ be returning   vectors. Their orbits are
relatively compact, since they are bounded and lie in
finite-dimensional subspaces $L(a_i)$. One can find a sequence
$n_k\to\infty$ such that $T^{n_k}a_i$ converges for each $i=1,\ldots
n$. Then for the sequence $m_k\to\infty$ of the form $n_{k+l}-n_k$
we have $T^{m_k}(a_i)\to a_i$. Thus if
$a=\lambda_1a_1+\ldots+\lambda_na_n$, then $T^{m_k}a\to a$ and $a$
is a returning vector.

So, the linear span $L$ of the set of returning   vectors itself
consists of returning   vectors. According to Lemma 1 $T:L\to L$ is
isometry. According to Theorem 1 $\dim L<\infty$. So, for each $x\in
X$ there exists $a\in L$ such that $T^nx\to O(a)\subseteq L$.

For every $x\in X$ the continuous function $\rho_x:L\to\Bbb R$
defined by the formula $\rho_x(a)=\liminf_n\|T^nx-T^na\|$ attains
its minimum $0$ at a unique point $a(x)\in L$. Clearly,
$\|T^nx-T^na(x)\|\to 0$. Linearity and boundedness $A:x\mapsto a(x)$
are obvious. Put $X_0=\ker A\subseteq X$, i.e. $x\in X_0
\Leftrightarrow T^nx\to 0$. The decomposition $X=X_0\oplus L$
corresponds to the condition (2). Theorem 2 is proved.

\vspace{1mm} {\it Proof of Theorem 3}: \ The set
$\tilde{K}=\cup_{t\in [0, 1]}T_t(K)\subseteq X$ is compact, since it
is the image of the compact set $K\times[0,1]$ under the map
$f(x,t)=T_tx$. It is easy to see that $\tilde{K}$ is occasionally
attracting for semigroup of powers  $\{T_1, T_2,\ldots\}$, i.e. the
operator $T_1:X\to X$ satisfies the conditions of Theorem 2. Let $L$
be a finite-dimensional subspace attracting $X$ under the action of
integer powers of $T_1$, i.e. $T_nx\to L$ for each $x\in X$. Show
that $T_tx\to L$ for every~$x$.

Suppose the contrary. Then there exist a number $\varepsilon>0$ and
a sequence $t_n\in \Bbb R, t_n\to\infty$ such that $\rho(T_{t_n}x,
L)>\varepsilon.$ Denote by $[t_n]$ and $\{t_n\}$ the integer and
fractional parts of the number $t_n$. It is possible to assume that
$\{t_n\}\to \beta\in[0,1]$ by switching to a subsequence. Then
 $$T_{t_n}x=T_{[t_n]+\{t_n\}}x=T_{[t_n]}T_{\{t_n\}}x\sim T_{[t_n]}T_\beta x\to L.$$
A contradiction. Theorem 3 is proved.\vskip1mm

\maketitle \centerline{\bf 3. Application to supercyclic
operators}\vskip1mm

%{\bf Lemma 4}. Let $T:X\to X$ be power bounded operator.
%If $T^na\not{\!\!\to}0$ and there
%exist $\lambda_k$ and $n_k$ such that $\lambda_kT^{n_k}a\to a$
%(in particular, if $a$ is supercyclic), then $a$ is a returning   vector.

{\it Proof of lemma 4}: \ If $T^na\not{\!\!\to}0$, then there exist
a \it bounded \rm sequence of scalars $c_k$ and a sequence of powers
$l_k\to\infty$ such that $c_kT^{l_k}a\to a$. Choose a subsequence
$m_k$ such that $c_k\to c$ and $cT^{m_k}a\to a$. Clearly, $|c|=1$.
In this case, $c^2T^{2m_k}a\to a$, $c^3T^{3 m_k}\to a$,... But 1 is
a limit point of the set $\{c^m\mid m\in\Bbb N\}$, therefore $a$ is
a limit point of the set  $\{T^{m\cdot n_k}a\mid m,k\in\Bbb
N\}$.$\Box$

{\it Proof of Theorem 4}: \ Rescaling $X$ by the norm (4), we may
suppose that $\|T\|\leq 1$. Let $a$ be supercyclic. In particular, a
is cyclic, i.e. the ${span}(O(a))$ is dense in X.

Assume that $T^na\not{\!\!\to}0$. By Lemma 4 $a$ is a returning
vector, therefore $T:X\to X$ is an~isometry by Lemma~1.

For any $x\in B_X$ there exist $\lambda_k$, $|\lambda_k|\leq 1$ and
$n_k\to\infty$ such that $\|\lambda_k T^{n_k}y- x\|\to 0$ or,
equivalent, $\|\lambda_ky-T^{-n_k} x\|\to 0$, therefore the
%one-dimensional
set $K=\{\lambda y\mid |\lambda|\leq 1\}$ is an occasionally compact
set for the \it isometry \rm $T^{-1}$. A contradiction with theorem
1.

Thus  $T^na\to 0$. But in this case $T^nx\to 0$ for every~$x$.
Indeed, for each $\varepsilon>0$, there is a vector of the form $c
T^k(a)$ that is $\varepsilon$-closed to $x$. Iterating $T$, we infer
$c T^{k+n}(a)\to_{n\to\infty} 0$; consequently, $\|T^n
x\|<\varepsilon$ for large $n$. Hence, $T^n x\to 0$. $\Box$\vskip1mm

I am grateful to Eduard Emelyanov, who attracted my attention to
papers \cite{LLY,Ph,S} and formulated a~hypothesis that led to
Theorems~2 and~3.

\end{document}